\documentclass[A4paper,11pt]{article}
\usepackage{amsfonts}
\usepackage{amsmath}
\usepackage{amssymb}
\usepackage{graphicx}
\usepackage{rlepsf}
\begin{document}
\title{Categorification of the Dichromatic Polynomial for Graphs}
\author{Marko Sto\v si\'c 
\thanks{The author is supported by {\it Funda\c c\~ao de Ci\^encia e
Tecnologia}/(FCT), grant no. SFRH/BD/6783/2001}\\
Departamento de Matem\'atica  and \\
CEMAT -  Centro de Matem\'atica e Aplica\c c\~oes\\
Instituto Superior T\'ecnico\\
Av. Rovisco Pais 1\\
1049-001 Lisbon\\ 
Portugal\\
e-mail: mstosic@math.ist.utl.pt
}

\date{}

\newtheorem{theorem}{Theorem}
\newtheorem{acknowledgment}[theorem]{Acknowledgment}
\newtheorem{algorithm}[theorem]{Algorithm}
\newtheorem{axiom}[theorem]{axiom}
\newtheorem{case}[theorem]{Case}
\newtheorem{claim}[theorem]{Claim}
\newtheorem{conclusion}[theorem]{Conclusion}
\newtheorem{condition}[theorem]{Condition}
\newtheorem{conjecture}[theorem]{Conjecture}
\newtheorem{corollary}[theorem]{Corollary}
\newtheorem{criterion}[theorem]{Criterion}
\newtheorem{definition}{Definition}
\newtheorem{example}[theorem]{Example}
\newtheorem{exercise}[theorem]{Exercise}
\newtheorem{lemma}{\indent Lemma}
\newtheorem{notation}[theorem]{Notation}
\newtheorem{problem}[theorem]{Problem}
\newtheorem{proposition}{Proposition}
\newtheorem{remark}[theorem]{Remark}
\newtheorem{solution}[theorem]{Solution}
\newtheorem{summary}[theorem]{Summary}
\newcommand{\ud}{\mathrm{d}}

\def\gcd{\mathop{\rm gcd}}
\def\Ker{\mathop{\rm Ker}}
\def\max{\mathop{\rm max}}
\def\map{\mathop{\rm map}}
\def\lcm{\mathop{\rm lcm}}
\def\kraj{\hfill\rule{6pt}{6pt}}
\def\diag{\mathop{\rm diag}}
\def\span{\mathop{\rm span}}
\def\deg{\mathop{\rm deg}}
\def\rank{\mathop{\rm rank}}
\def\sgn{\mathop{\rm sgn}}
\def\kvn{\{n\}}
\def\F{\mathbb{F}}
\def\R{\mathbb{R}}
\def\C{\mathcal{C}}
\def\K{\mathbb{K}}
\def\Z{\mathbb{Z}}
\def\Q{\mathbb{Q}}
\arraycolsep 6pt

\maketitle

\begin{abstract}
For each graph and each positive integer $n$, we define a chain complex
whose graded Euler characteristic is equal to an appropriate $n$-specialization
of the dichromatic polynomial. This also gives a categorification of
$n$-specializations of the Tutte polynomial of graphs. Also, for each graph and integer $n\le 2$, we define the different one variable $n$-specializations of the dichromatic polynomials, and for each polynomial we define graded chain complex whose graded Euler characteristic is equal to that polynomial. Furthermore, we explicitly categorify the specialization of the Tutte polynomial for graphs which corresponds to the Jones polynomial of the appropriate alternating link. 
\end{abstract}

\section{Introduction}
One of the most interesting and promising recent developments in knot theory, 
and in mathematics in general,
is the ``categorification" of link invariants, initiated by Khovanov in \cite{kov}.
For each link $L$ in $S^3$ he defined a graded chain complex, with grading 
preserving differentials, whose graded Euler characteristic is equal to
the Jones polynomial of the link $L$. This is done by starting from the state sum 
expression for the Jones polynomial (which is written as an alternating sum),
then constructing for each term a module whose graded dimension is 
equal to the value of that term, and finally, constructing the
differentials 
as appropriate grading preserving maps, so that 
the complex obtained is a link invariant. There
are similar constructions for the categorification of various 
link polynomial invariants (\cite{kovroz}, \cite{kov3}, \cite{kol}, \cite{bn1}).\\
\indent Recently, a similar construction has also been done for the chromatic polynomial
for graphs \cite{gr}. In this paper we will generalize their work, by categorifying
the dichromatic polynomial (and consequently the Tutte polynomial). Since
the dichromatic polynomial is a two-variable polynomial and the standard techniques
of categorification work with one-variable polynomials, we will define a
class of one-variable specializations $P_n(G)$, one for every integer $n\le 2$ - 
a family which is enough to recover the original dichromatic
polynomial - and
then for each of them we construct a chain complex whose graded Euler characteristic
is equal to $P_n(G)$. \\
\indent As is well-known, to each alternating link $L$ there
corresponds bijectively a planar graph
 $G(L)$, 
and the value of the Jones polynomial of the link $L$ corresponds to
the specialization of the
Tutte polynomial $T(x,1/x)$ or analogously to the specialization of the dichromatic polynomial
$J(q)=P(q,q^2/(q-1))$ (see e.g. \cite{bol} or \cite{kauf}). This is exactly one of the specializations from the
 previous paragraph (the one for $n=2$), and hence we obtain the categorification of the polynomial $J(q)$. However, we also give the alternative description in terms of the enhanced states. \\
\indent Even more, we define different set of one-variable specializations $Q_n(G)$, one for every positive integer $n$, and then categorify each of these one-variable polynomials. Although this set of specializations ``miss" the Jones polynomial, it has the advantage that the chain groups of the complexes are finite-dimensional.\\
\indent The organization of the paper is as follows: in Section 2 we recall the definitions 
of the dichromatic polynomial, define the appropriate specializations
and explain how we will define the chain complexes. In Section 3, we explicitly define
the chain complex for the first specialization. Further, in Section 4, we construct the chain complex whose 
graded Euler characteristic is equal to $J(q)$, the direct correspondent of the Jones polynomial
for alternating links, in terms of the enhanced states. Finally, in Section 5, we categorify the second set of one-variable specializations.

\section{Preliminaries}

Let  $G$ be a graph specified by a set of vertices $V(G)$ and a set of
edges $E(G)$. If $e\in E(G)$ is an arbitrary edge of the graph $G$,
then by $G-e$ we denote the graph $G$ with the edge $e$ deleted, and
by $G/e$ the graph obtained by contracting edge $e$ (i.e. by
identifying the vertices incident to $e$ and deleting $e$). The
dichromatic polynomial of the graph, $P_G(q,v)$, is the two-variable
generalization of the
chromatic polynomial of the graph. It is given by the following axioms:
\begin{eqnarray*}
(A1)&\quad P_G=P_{G-e}-qP_{G/e},\\
(A2)&\quad P_{N_k}=v^k,
\end{eqnarray*}
where $N_k$ is the graph with $k$ vertices and no edges.\\
\indent These two axioms determine the polynomial uniquely. Obviously, if we put $q=1$, we obtain the usual chromatic polynomial. Furthermore, from $(A1)$ we have a recursive expression for the dichromatic polynomial in terms of the value of the polynomial on graphs with a smaller number of edges. By repeated use of $(A1)$ we will obtain the value of the dichromatic polynomial as a sum of  contributions from all spanning subgraphs of $G$ (subgraphs that contain all vertices of $G$), which we will call states.
Furthermore, if for each subset $s\subset E(G)$ we denote by $[G:s]$ the graph whose set of vertices is $V(G)$ and set of edges is $s$, then the contribution of the graph $[G:s]$ is $(-1)^{|s|}q^{|s|}v^{k(s)}$, where $|s|$ is the number of elements of $s$ and $k(s)$ is the number of connected components of $[G:s]$. Hence, we obtain the expression:
$$P_G(q,v)=\sum_{s\subset E(G)}{(-1)^{|s|}q^{|s|}v^{k(s)}}=\sum_{i\ge 0}{(-1)^i q^i\sum_{s\subset E(G),|s|=i}{v^{k(s)}}},$$
which is called the state sum expansion of the polynomial $P_G$.\\
\indent For the Tutte polynomial $T_G(x,y)$ of a graph $G$ there exists a similar state-sum expression given by:
$$T_G(x,y)=\sum_{s\subset E(G)}{(x-1)^{k(s)-k(E(G))}(y-1)^{|s|-N+k(s)}},$$
where by $N$ we denote the number of vertices of $G$. It is obviously related to the dichromatic polynomial by:
\begin{eqnarray*}
T_G(x,y)&=&\sum_{s\subset E(G)}{(x-1)^{k(s)-k(E(G))}(y-1)^{|s|-N+k(s)}}=\\
        &=&(x-1)^{-k(E(G))}(y-1)^{-N}\sum_{s\subset E(G)}{(y-1)^{|s|}{(x-1)(y-1)}^{k(s)}}=\\
        &=&(x-1)^{-k(E(G))}(y-1)^{-N}P_G(1-y,(x-1)(y-1)).
\end{eqnarray*}
Hence, we have that the Tutte polynomial $T_G(x,y)$ is a multiple of the dichromatic polynomial $P_G(q,v)$, when we take $q=1-y$ and $v=(x-1)(y-1)$.\\ 

\indent Our aim is to define a graded chain complex whose graded Euler characteristic is equal to the dichromatic polynomial. However, since the dichromatic polynomial is a two-variable polynomial, we will first define an infinite set of one-variable specializations and then ``categorify" each of the specializations. This situation is very similar to the two-variable HOMFLY polynomial for knots and its infinite set of one-variable specializations (one for each positive integer $n$, with $n=2$ being the Jones polynomial). Note that in this way we will also categorify the Tutte polynomial.\\

\indent In order to do this we will introduce new variables. We will assume that $q > 1$, and introduce a new variable $a$ as $a=v(q-1)$, i.e. $v=a/(q-1)$. In this way, we obtain a two-variable polynomial $\tilde{P}_G(q,a)$, and the one-variable specializations we will define to be $P_{G,n}(q)=\tilde{P}_{G}(q,q^n)$, for every integer $n \le 2$. Note that we then have $v=q^n/(q-1)=q^{n-1}/(1-q^{-1})$. We will denote the expression $q^n/(q-1)=q^{n-1} \sum_{i \ge 0}{q^{-i}}$ by $\{n\}$.\\

\begin{figure}[h]
\centerline{\relabelbox
\epsfysize 10cm
\epsfbox{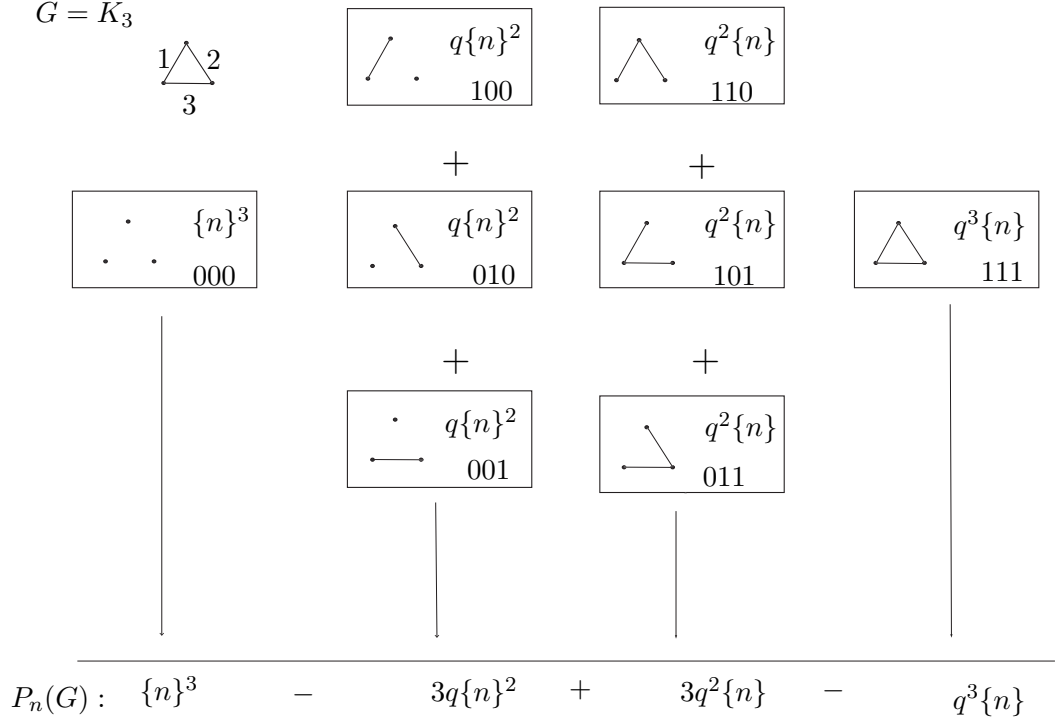}
\adjustrelabel <-15pt,20pt> {G=K3}{$G=K_3$}
\adjustrelabel <-4pt,0pt> {1}{$1$}
\relabel {2}{$2$}
\adjustrelabel <-3pt,-6pt> {3}{$3$}
\relabel {000}{$000$}
\relabel {001}{$001$}
\relabel {010}{$010$}
\relabel {100}{$100$}
\relabel {011}{$011$}
\relabel {101}{$101$}
\relabel {110}{$110$}
\relabel {111}{$111$}
\adjustrelabel <-5pt,-7pt> {n3}{$\{n\}^3$}
\adjustrelabel <-10pt,-7pt> {qn1}{$q\{n\}^2$}
\adjustrelabel <-10pt,-7pt> {qn2}{$q\{n\}^2$}
\adjustrelabel <-10pt,-7pt> {qn3}{$q\{n\}^2$}
\adjustrelabel <-5pt,-7pt> {q1n}{$q^2\{n\}$}
\adjustrelabel <-6pt,-7pt> {q2n}{$q^2\{n\}$}
\adjustrelabel <0pt,-7pt> {q3n}{$q^2\{n\}$}
\adjustrelabel <-10pt,-7pt> {q3n1}{$q^3\{n\}$}
\relabel {+1}{{\Large $+$}}
\relabel {+2}{{\Large $+$}}
\relabel {+3}{{\Large $+$}}
\relabel {+4}{{\Large $+$}}
\relabel {pnG:}{$P_n(G):$}
\relabel {v1}{$\{n\}^3$}
\relabel {v2}{$3q\{n\}^2$}
\relabel {v3}{$3q^2\{n\}$}
\relabel {v4}{$q^3\{n\}$}
\relabel {-}{$-$}
\relabel {+}{$+$}
\relabel {-1}{$-$}
\endrelabelbox}
\caption{State sum expression for the $n$-specialization of the dichromatic polynomial}
\label{sl1}
\end{figure}

\indent For every integer $n \le 2$ and graph $G$, we will define a chain complex whose graded Euler characteristic is equal to $P_{G,n}$. Denote by $m$ the number of edges of the graph $G$. The construction will depend on the ordering of the edges of $G$ ($e_1,\ldots,e_m$). Then each subset $s_{\epsilon}\subset E(G)$ will be uniquely determined by an element $\epsilon=(\epsilon_1,\ldots,\epsilon_m)\in \{0,1\}^m$, 
where $\epsilon_i=1$ if $e_i \in s$ and if $\epsilon_i=0$ if $e_i \notin s$. Obviously, each such $m$-tuple is uniquely determined by the subset $s$ of $E(G)$, and so we have a bijective correspondence between $\{0,1\}^m$ and the set of all $s\subset E(G)$. Thus, every element $\epsilon \in \{0,1\}^m$ determines a set $s_{\epsilon}$ of edges of $G$, and hence a graph $[G:s_{\epsilon}]$ which we will denote by $G_{\epsilon}$.  \\
\indent We will construct our complex by ``flattening" the cube of states (see figure \ref{sl1}).

To each vertex $\epsilon$ of the cube $\{0,1\}^m$ (which is skewed 
such that in the $i$-th column, $0\le i \le m$, are all the $\epsilon$'s 
such that $|\epsilon|=i$), we will assign a graded $\Z$-module 
whose graded dimension is equal to the contribution of 
the subgraph $G_{\epsilon}$. In other words, a module whose 
graded dimension is $q^{|\epsilon|}\kvn^{k_{\epsilon}}$, where
$k_{\epsilon}$ is the number of connected components of $G_{\epsilon}$.
We will build chain groups by taking direct sums along the columns.
Finally, we will define differentials such that they are grading preserving,
and in this way we will obtain a chain complex whose graded Euler characteristic 
is equal to $n$-th specialization $P_{G,n}$ of the dichromatic polynomial.\\
\indent In the next section we will define the chain complex more precisely.   \\

\indent Also, in Section 5, we will observe the different set of specializations of the dichromatic polynomial. Namely, for each positive integer $n$ we define:
$$Q_n=P(q^n,1+q+\ldots+q^n).$$
We also obtain analogous state-sum expression as for $P_n$, and then define the cubic complex as outlined before.

\section{The cubic complex construction of the chain complex}

\indent We are going to assign to each state a graded $\Z$-module whose graded dimension is $\kvn^{k}$. Before going to the definition, we will recall some basic facts about graded $\Z$-modules with integer gradings and their graded dimensions.
\subsection{Graded dimension of a graded $\Z$-module}
\begin{definition}
Let $M=\oplus_{k}M_k$ for $k \in \Z$ be a graded $\Z$-module where $\{M_k\}$ denotes the $k$th graded component of $M$. The graded dimension of $M$ is the power series
$$q\dim M:=\sum_k{q^k\rank(M_k)},$$
where $\rank(M_k)=\dim_{\Q}(M_k\otimes\Q)$.
\end{definition}

\indent The direct sum and tensor product can be defined in the graded category in an obvious way. The following proposition is straightforward.
\begin{proposition}
Let $M$ and $N$ be graded $\Z$ modules. Then
\begin{eqnarray*}
q\dim (M\oplus N)&=&q\dim(M)+q\dim(N) \\
q\dim(M\otimes N)&=&q\dim(M) q\dim(N).
\end{eqnarray*}
\end{proposition}
\begin{definition}
Let $\{\l\}$, $l \in \Z$, be the ``degree shift" operation on graded $\Z$-modules. That is, if $M=\oplus_k M_k$ is a graded $\Z$ module where $M_k$ denotes the $k$th graded component of $M$, we set $M\{l\}_k:=M_{k-l}$. Then we have 
$q\dim M\{l\}=q^{l}q\dim M.$ In other words, all the degrees are increased by $l$.
\end{definition}

\begin{example}
Let $M=\Z[x_i]\{n-1\}$ be the ring of polynomials in one variable. If we define the degree of the variable $x_i$ to be $-1$, then $M$ becomes graded free $\Z$-module, with quantum dimension $q\dim M =q^{n-1}\sum_{i\ge 0}{q^{-i}}=\kvn$. Further, if $k>0$, then $M^{\otimes k}$ is isomorphic to the ring of polynomials in $k$ variables $M_k=\Z[x_1,\ldots,x_k]\{k(n-1)\}$, and we have that $q\dim M^{\otimes k} = q\dim M_k = {\{ n\}}^k$.  
\end{example}
\begin{example}
Let $V$ be the graded free $\Z$-module with $n+1$ basis elements: $X^0(=1),X,X^2,\ldots,X^n$ such that the degree of $X^i$ is equal to $n-i$, for $i=0,\ldots,n$. Then we have $q\dim M=1+q+\ldots+q^n$. Furthermore, we have that $V^{\otimes k}\{l\}=q^l(1+q+\ldots+q^n)^k.$
Notice that we can also describe $V$ as the quotient of the ring of polynomials by $V=\Z[X]/(X^{n+1})$.
\end{example}

\indent

\indent We are now ready to explain our construction. Let $G$ be a graph with $m$ ordered edges and let $M$ be as in Example 1. To each vertex $\epsilon=(\epsilon_1,\ldots,\epsilon_m)$ of the cube $\{0,1\}^m$, we associate the graded free $\Z$-module $M_{\epsilon}$ in the following way: First of all, we order the vertices of $G$, and to the $i$-th one assign the index $i$. Let $K_1$, $K_2$,$\ldots$,$K_l$ be the set of connected components of $G_{\epsilon}$. Then to the component $K_j$ we will assign variable $x_{i_j}$, where $i_j$ is the smallest index among the vertices $v$ that belong to $K_j$. Finally, we define $M_{\epsilon}=\Z[x_{i_1},\ldots,x_{i_l}]\{l(n-1)+|\epsilon|\}$.\\In other words, to each connected component of $G_{\epsilon}$ we assign a copy of $M=\Z[x]\{n-1\}$, then take the tensor product and finally increase the degree by $|\epsilon|$. From the definition we have that $q\dim M_{\epsilon}$ is (up to the sign $(-1)^{|\epsilon|}$) the contribution of $P_{G_{\epsilon}}$ to the dichromatic polynomial.\\

\begin{figure}[h]
\centerline
\relabelbox
\epsfysize 10cm
\epsfbox{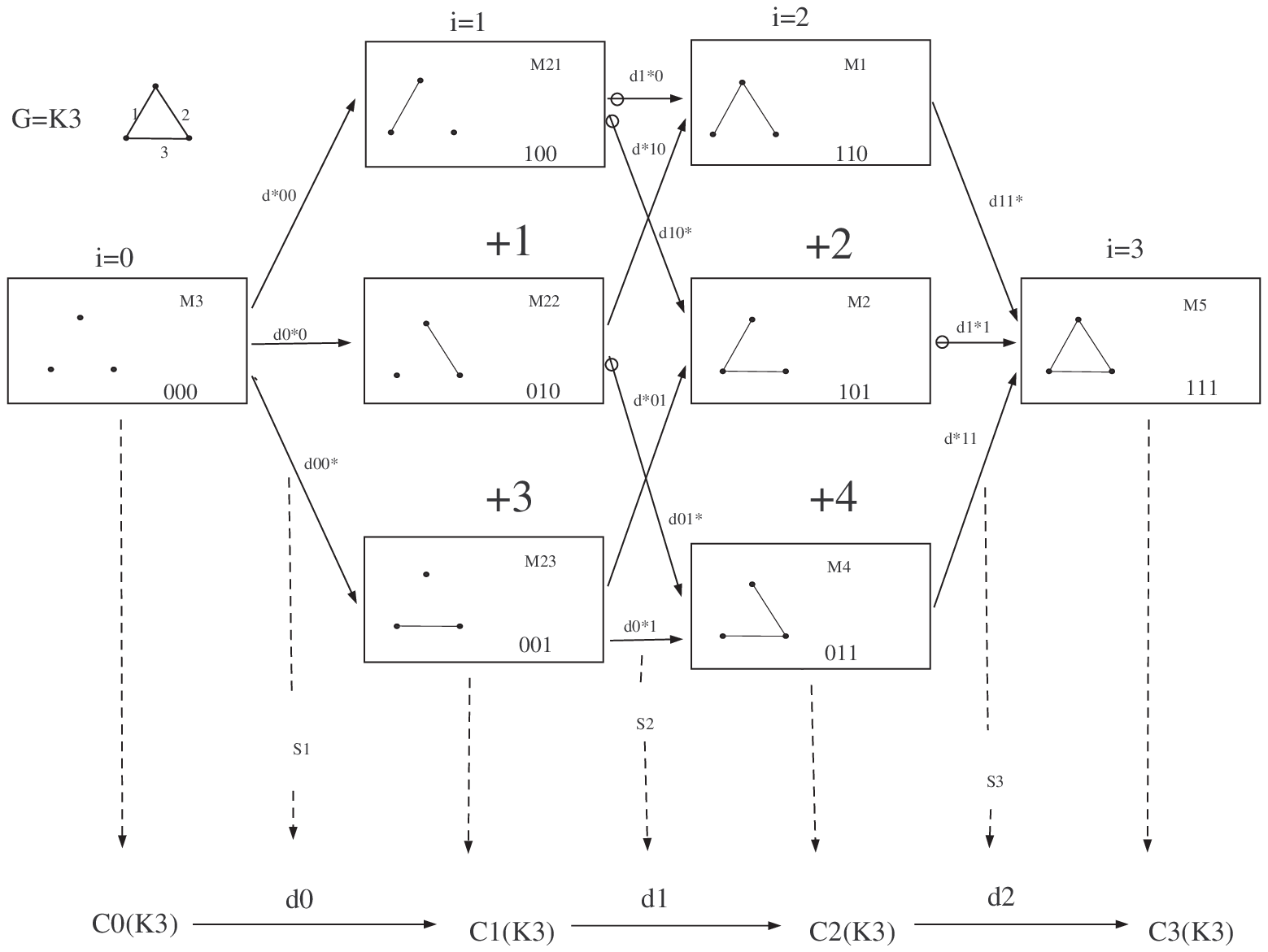}
\adjustrelabel <10pt,20pt> {G=K3}{$G=K_3$}
\adjustrelabel <-6pt,0pt> {1}{$1$}
\relabel {2}{$2$}
\adjustrelabel <-3pt,-5pt> {3}{$3$}
\relabel {000}{$000$}
\relabel {001}{$001$}
\relabel {010}{$010$}
\relabel {100}{$100$}
\relabel {011}{$011$}
\relabel {101}{$101$}
\relabel {110}{$110$}
\relabel {111}{$111$}
\adjustrelabel <-7pt,-7pt> {M3}{$M^{\otimes 3}$}
\adjustrelabel <-20pt,-7pt> {M21}{$M^{\otimes 2}\{1\}$}
\adjustrelabel <-20pt,-7pt> {M22}{$M^{\otimes 2}\{1\}$}
\adjustrelabel <-20pt,-7pt> {M23}{$M^{\otimes 2}\{1\}$}
\adjustrelabel <-5pt,-6pt> {M1}{$M\{2\}$}
\adjustrelabel <-6pt,-7pt> {M2}{$M\{2\}$}
\adjustrelabel <-1pt,-7pt> {M4}{$M\{2\}$}
\adjustrelabel <-10pt,-7pt> {M5}{$M\{3\}$}
\relabel {+1}{$\bigoplus$}
\relabel {+2}{$\bigoplus$}
\relabel {+3}{$\bigoplus$}
\relabel {+4}{$\bigoplus$}
\relabel {i=0}{$i=0$}
\relabel {i=1}{$i=1$}
\relabel {i=2}{$i=2$}
\relabel {i=3}{$i=3$}
\adjustrelabel <-8pt,0pt> {C0(K3)}{$C^0(K_3)$}
\adjustrelabel <-6pt,0pt> {C1(K3)}{$C^1(K_3)$}
\adjustrelabel <-6pt,0pt> {C2(K3)}{$C^2(K_3)$}
\adjustrelabel <-5pt,0pt> {C3(K3)}{$C^3(K_3)$}
\relabel {d0}{$d_0$}
\relabel {d1}{$d_1$}
\relabel {d2}{$d_2$}
\adjustrelabel <-5pt,0pt> {d*00}{$d_{\ast 00}$}
\adjustrelabel <0pt,2pt> {d0*0}{$d_{0\ast 0}$}
\relabel {d00*}{$d_{00\ast}$}
\relabel {d1*0}{$d_{1\ast 0}$}
\relabel {d*10}{$d_{\ast 10}$}
\relabel {d10*}{$d_{10\ast}$}
\relabel {d*01}{$d_{\ast 01}$}
\relabel {d01*}{$d_{01\ast}$}
\relabel {d1*1}{$d_{1\ast 1}$}
\relabel {d0*1}{$d_{0\ast 1}$}
\relabel {d11*}{$d_{11\ast}$}
\adjustrelabel <-3pt,0pt> {d*11}{$d_{\ast 11}$}
\adjustrelabel <-30pt,0pt> {S1}{$\sum_{\xi=0}{(-1)^{\xi}d_{\xi}}$}
\adjustrelabel <-30pt,0pt> {S2}{$\sum_{\xi=1}{(-1)^{\xi}d_{\xi}}$}
\adjustrelabel <-30pt,0pt> {S3}{$\sum_{\xi=2}{(-1)^{\xi}d_{\xi}}$}
\endrelabelbox
\caption{The chain complex and the differentials}
\label{sl2}
\end{figure}

\indent To get the chain groups we ``flatten" the cube by taking direct sums along the columns. More precisely:
\begin{definition}
We set the $i$th chain group $C_n^i(G)$ of the chain complex $\C_n(G)$ to be the direct sum of all $\Z$-modules at height $i$, i.e. $C_n^i(G)=\oplus_{|\epsilon|=i}{M_{\epsilon}(G)}$.\\
The grading is given by the degree of the elements and we can write the $i$th chain group as $C_n^i(G)=\sum_{j\in\Z}C^{i,j}_n(G)$, where $C^{i,j}_n(G)$ denotes the set of elements of degree $j$ of $C_n^i(G)$.
\end{definition}

\subsection{Graded chain complex, graded Euler characteristic}

\begin{definition}
Let  $M=\oplus_j M_j$ and $N=\oplus_j N_j$ be graded $\Z$-modules where $M_j$ and $N_j$ denote the $j$th graded component of $M$ and $N$, respectively. A $\Z$-module map $\alpha:M\rightarrow N$ is said to be graded with degree $d$ if $\alpha(M_j)\subset N_{j+d}$, i.e. elements of degree $j$ are mapped to elements of degree $j+d$. A $\Z$-module map is called degree preserving if it is graded of degree zero.\\
A graded chain complex is a chain complex for which the chain groups are graded $\Z$-modules and the differentials are graded.
\end{definition}

\begin{definition}
The graded Euler characteristic $\chi_q(\C)$ of a graded chain complex $\C$ is the alternating sum of the graded dimensions of its homology groups, i.e. $\chi_q(C)=\sum_{0\le i\le m}{(-1)^i q\dim(H^i)}$.
\end{definition}

\begin{proposition}
(\cite{bn}, \cite{gr}) If the differential is degree preserving and all $\C^{i,j}$'s (sets of fixed of degree of each chain group) are finite dimensional, the graded Euler characteristic is also equal to the alternating sum of the graded dimensions of its chain groups, i.e.
$$\chi_q(C)=\sum_{0\le i\le m} (-1)^iq\dim(H^i)=\sum_{0\le i\le m}(-1)^iq\dim (C^i).$$
\end{proposition}

\indent From the last Proposition we have that if we define a degree preserving differential between the groups $C_n^i(G)$, the chain omplex $\C_n(G)$ obtained in this way, will have as Euler characteristic, the specialization $P_{G,n}(q)$ of the dichromatic polynomial. Hence, we are left with defining such differential.

\subsection{The differential}

\indent We define the differential in the, now, standard way. We are first going to define per-edge maps between some vertices of the cube $\{0,1\}^m$ - the maps that correspond to the edges of the cube. We define them as  linear maps and such that the cube is commutative, i.e.  every square is commutative. Then we build the differential by summing with appropriate signs along columns, and hence obtain a map whose square is zero. This is presented graphically in figure \ref{sl2}.\\
\indent So, let us first define per-edge maps. Each vertex of the cube $\{0,1\}^m$ is labeled with some $\epsilon=(\epsilon_1,\ldots,\epsilon_{m})\in\{0,1\}^m$. There are maps between two vertices only if one of the markers $\epsilon_i$ is changed from $0$ to $1$ when one goes from the first vertex to the second vertex and all the other $\epsilon_i$ are unchanged.\\
\indent Denote by $\epsilon$ the label of the first vertex. If the marker which is changed from $0$ to $1$ has index $j$ then the map will be labeled $d_{\epsilon'}$, where $\epsilon'=(\epsilon'_1,\ldots,\epsilon'_{m})$ with $\epsilon'_i=\epsilon_i$ if $i\ne j$ and $\epsilon'_i=\ast$ if  $i=j$. Denote by $l$ the sum of components of $\epsilon$, and by $l$ the number of connected components of $G_{\epsilon}$.
Changing exactly one marker from 0 to 1 corresponds to adding an edge. \\
\indent If adding that edge does not affect the number of components, then it means that the added edge will belong to some component, say the $p$-th one, and that the remaining components will remain unchanged. In this case, we have to define the grading preserving map from $M_{\epsilon}=\Z[x_{i_1},\ldots,x_{i_l}]\{l(n-1)+|\epsilon|\}$ to $M_{\epsilon'}=\Z[x_{i_1},\ldots,x_{i_l}]\{l(n-1)+|\epsilon|+1\}$. We define the map as the multiplication by $x_{i_p}$. \\
\indent If adding the edge decreases the number of components by one, then we have that two components (say $p$-th and $q$-th one) are merging into one, and the remaining $(l-2)$ components are left the same. Suppose that $i_p<i_q$. We have to define a degree preserving map from $M_{\epsilon}=\Z[x_{i_1},\ldots,x_{i_l}]\{l(n-1)+|\epsilon|\}$ to $M_{\epsilon'}=\Z[x_{i_1},\ldots,\hat{x}_{i_q},\ldots,x_{i_l}]\{(l-1)(n-1)+|\epsilon|+1\}$, where by $\hat{x}_{i_q}$ we have indicated that $x_{i_q}$ is omitted. We define it on the basis monomials by:
$$\prod_{j=1,\ldots,l}{x_{i_j}^{a_j}} \mapsto \prod_{{j=1,\ldots,l} \atop {j\ne q}}{x_{i_j}^{b_j}},$$
where $b_i=a_i$, for $i\ne p$, and $b_p=a_p+a_q+2-n$. In other words, we map $x_{i_q}$ to $x_{i_p}$ and multiply the result by $x_{i_p}^{2-n}$.
 
\indent The per-edge differentials obtained in this way obviously make the cube commutative. We will sprinkle signs around on the edges of the cube and thus obtain a map whose square is zero. Namely, we define 
the differential $d^i:C_n^i(G)\rightarrow C^{i+1}_n(G)$ on the chain complex $\C_n(G)$ by $$d^i:=\sum_{|\epsilon|=i}{(-1)^{\epsilon}d_{\epsilon}},$$ 
where we sum over all $m$-tuples of $0$'s, $1$'s, and exactly one $\ast$. Here we denote the number of $1$'s in $\epsilon$ by $|\epsilon|$, and $(-1)^{\epsilon}$ equals $-1$ if there are an odd number of $1$'s before $\ast$ in $\epsilon$, and $+1$ if there are an even number of them. In Figure \ref{sl2} we indicated the differentials with minus sign by adding a small circle at the tail of the arrow.\\

\indent A straightforward calculation implies:

\begin{proposition}
This defines a differential, that is, $d^2=0$.
\end{proposition}

\indent Now we really have a chain complex $\C_n(G)$ where the chain groups and the differential are defined as in the previous two paragraphs, and according to Proposition 2, we have

\begin{theorem}
The Euler characteristic of this chain complex $\C_n(G)$ is equal to the $n$-specialization $P_{G,n}$ of the dichromatic polynomial of the graph $G$.
\label{t1}
\end{theorem}

\indent Even though our construction depends on the ordering of the edges of the graph, in exactly the same way as in \cite{gr}, section 2.2.3, we obtain the following:

\begin{theorem}
The homology groups of the chain complex $\C_n(G)$ are graph invariants.
\label{t2}
\end{theorem}

\indent For each graph $G$ and integer $n\le 2$, we define the (2-variable) Poincar\'e polynomial of the corresponding chain complex
$\C_n(G)$, in the standard way, by $R_{G,n}(q,t)=\sum_{i\in \Z}{t^i q\dim (H^i)}$. Then from Theorems \ref{t1} and 
\ref{t2} we obtain:

\begin{proposition}
a) For every integer $n\le 2$, the polynomial $R_{G,n}$ is a graph invariant.\\
b) The $n$-specialization of the dichromatic polynomial $P_{G,n}(q)$ is equal to $R_{G,n}(q,-1)$. 
\end{proposition}
 
\section{Categorification of the Jones polynomial for alternating links}

As is well-known, every planar graph is in bijective correspondence with some knot universe (knot shadow), which in turn gives rise to an alternating link diagram. As is also known (see e.g. \cite{bol}), the Jones polynomial of an alternating link is related (equal up to a multiple) with the value of the Tutte polynomial $T_G(x,1/x)$ of a planar graph corresponding to it. As we saw in the introduction the value of this specialization of Tutte polynomial is (up to a multiple) equal to $P_G(q,q^2/(q-1))$ (with $q=1-1/x$). This can also be seen directly (see e.g. \cite{kauf}). Hence, we can obtain a categorification of the Jones polynomial of an alternating link, by categorifying the one variable specialization $J_G(q)=P_G(q,q^2/(q-1))$ of the dichromatic polynomial.
As we saw in the previous section, we managed to categorify an infinite set of specializations of the dichromatic polynomial and for $n=2$ we obtain the categorification of $J_G(q)$. However, in this section we give alternative description of the chain complex in terms of enhanced states.\\
\indent First of all, notice that we can write the value $v$ (which in the case of $J_G(q)$ is equal to $q^2/(q-1)$) as the following series (for $|q|>1$):
\begin{equation}
v=q/(1-q^{-1})=q\sum_{i\ge 0}q^{-i}=\sum_{i\ge 0}q^{1-i},
\label{v}
\end{equation}
and with that value of $v$, the polynomial $J_G(q)$ is given by 
$$J_G(q)=\sum_{s\subset E(G)}{(-1)^{|s|}q^{|s|}v^{k(s)}}.$$ 
\indent In order to  define the graded complex, we will introduce the notion of enhanced states of a graph $G$. An enhanced state $S$ is given by a pair $S=(s,l)$, where $s$ is a state (subset of the set of edges $E(G)$), and $l$ is a labeling of the connected components of $[G:s]$ by nonnegative integers, i.e. a sequence of $k(s)$ nonnegative integers $l_1,\ldots,l_{k(s)}$. Denote by $|l|$, the sum of all $l_i$'s, $i=1,\ldots,k(s)$. To each enhanced state $S$ we will assign two numbers: $i(S)=|s|$, and $j(S)=|s|+k(s)-|l|$. Now, we define $C^{i,j}(G)$ ($j$-th graded component of the $i$-th chain group, $0\le i \le m$, $j\in \Z$) as the span over $\Z$ of all enhanced states $S$ of $G$ such that $i(S)=i$ and $j(S)=j$.\\
\indent We will define a (degree preserving) differential $d:C^{i,j}(G)\rightarrow C^{i+1,j}(G)$, by giving its value on each enhanced state $S=(s,l)\in C^{i,j}(G)$. We define:
\begin{equation}
d(S)=\sum_{e\in E(G)\setminus s}{(-1)^{n_e}}S_e,
\label{d}
\end{equation}
with $n_e$ being the number of edges in $s$ ordered before $e$, and $S_e=(s_e,l_e)$ is an enhanced state whose state is $s_e=s \cup \{e\}$ and whose labeling $l_e$ is defined as follows: denote by $E_1,\ldots,E_k$ the connected components of $[G:s]$.\\
\indent If $e$ connects some $E_i$ to itself then the components of $[G:s_e]$ are $E_1,\ldots,E_i\cup \{e\},\ldots,E_k$, and we define the labeling $l_e$ by $l_e(E_i\cup \{e\})=l(E_i)+1$ and $l_e(E_j)=l(E_j)$, for $j \ne i$.\\ 
\indent If $e$ connects some $E_i$ to $E_j$ (say $E_1$ and $E_2$), then the components of $[G:s_e]$ are $E_1\cup E_2 \cup \{e\}, E_3,\ldots, E_k$, and we define the labeling $l_e$ by $l_e(E_1 \cup E_2 \cup \{e\})=l(E_1)+l(E_2)$ and $l_e(E_j)=l(E_j)$, for $j>2$.\\

\begin{proposition}
The mapping $d$, defined by (\ref{d}), satisfies $d^2=0$.
\end{proposition}

\indent \textit{Proof:} \\
\indent Let $i$ and $j$ be arbitrary integers, and let $S$ be an enhanced state from $C^{i,j}$. If $e$ and $f$ are two distinct edges from $E(G)\setminus s$, then from the definition of the states $S_e$ in the image of the differential, we have $(S_e)_f=(S_f)_e$. Suppose that we have introduced an ordering $<$ on the set $E(G)$. Also, denote by $n'_e$, the number of edges in $s\cup \{f\}$ ordered before $e$, and by $n'_f$, the number of edges in $s\cup \{e\}$ ordered before $f$. Now we have:
\begin{eqnarray*}
d(d(S))&=&d(\sum_{e\in E(G)\setminus s}{(-1)^{n_e}}S_e)= \\
       &=&\sum_{f \in E(G)\setminus s_e}{(-1)^{n'_f}}{(-1)^{n_e}}(S_e)_f=\\
       &=&\sum_{{\scriptsize{\begin{array}{c}e,f \in E(G)\setminus s \\ e<f           \end{array}}}}{((-1)^{n'_f}(-1)^{n_e}+(-1)^{n'_e}{(-1)^{n_f}}})(S_e)_f=0,
\end{eqnarray*} 
as required. \kraj
\\
\indent In this way we have obtained the sequence of chain complexes (one for each degree $j$). We could take the direct sum along columns (fixed $i$), and obtain one graded chain complex $\C(G)$ with the chain groups given by:
$$C^i(G)=\oplus_{j\in \Z}C^{i,j}(G).$$
The graded dimension of $C^i(G)$ is equal to the sum of the graded dimensions corresponding to each state $s$ with $|s|=i$. From the definition of the gradings of our enhanced states we have that the graded dimension corresponding to each such state $s$ is equal to 
$$q^i{(\sum_{j\ge 0}{q^{1-j}})}^{k(s)}=q^iv^{k(s)},$$
with $v$ given by (\ref{v}).
Furthermore, since we have defined degree  preserving differentials, we have the following:
\begin{theorem}
The graded Euler characteristic of the chain complex $\C(G)$ is equal to $J_G(q)$.
\end{theorem}  

\section{Finite dimensional set of specializations}

In this section we define the alternative set of one-variable specializations and the corresponding set of chain complexes that categorify them. Although this set of specializations ``miss" the Jones polynomial from the previous section, the advantage is that the chain groups are finite dimensional.\\
\indent In this case we parametrize the set of specializations by positive integers. For each positive integer $n$, we define the polynomial $Q_{G,n}$ by
$$Q_{G,n}(q)=P_G(q^n,1+q+q^2+\ldots+q^n).$$
Note that we have $v=1+q+\ldots+q^n=(q^{n+1}-1)/(q-1)$. Further on, we will denote the expression $1+q+\ldots+q^n$ by $q_n$.\\
\indent The categorification of $Q_n$ follows the same lines as the categorification in Section 3. Namely, we organize the summand of the state sum expression for $Q_n$ in the same cube as in the figure \ref{sl1}, with the replacement of $\{n\}$ with $q_n$ and of $q$ with $q^n$. Furthermore, to each $\epsilon \in \{0,1\}^m$ (i.e. to each vertex of the cube of resolutions) we assign $\Z$-module $V_{\epsilon}(G)=V^{\otimes k(\epsilon)}\{n|\epsilon|\}$, with $V$ being the graded $\Z$-module defined in Example 2 in Section 3.1. In other words, we assign one copy of $V$ to each connected component of the resolution of the graph. In terms of pictures, we make the same picture as the figure \ref{sl2}, only with replacing $M$ with $V$ and all shifts $\{i\}$ by $\{ni\}$. Finally the chain groups $D^i_n(G)$ are obtained by summing along the columns, i.e. $D^i_n(G)=\oplus_{|\epsilon|=i}V_{\epsilon}(G)$.\\
\indent Now, we move to the differential. Prescription is the same as in the previous sections. We only have to define per-edge maps and make the cube commutative, and then the signs (defined in the same way as previously) will make the cube anti-commutative, which after summing over columns defines the map whose square is equal to zero, i.e. the differential.\\
\indent Again, the edges of the cube of resolutions correspond to adding of an edge of the graph G. This can affect the connected components of the resolutions in two ways: either the added edge connects two components (and consequently decreases the number of components by 1) or the added edge belongs to one of the components (and hence preserves the number of components).\\
\indent In the first case, we have that two components are merging into one (say $p$-th and $q$-th one) and the remaining $l-2$ components are left the same. We have to define a degree preserving map from $V^{\otimes l}\{ni\}$ to  $V^{\otimes (l-1)}\{n(i+1)\}$. We define it to be the identity on the $l-2$ tensor factors corresponding to the connected components which are left unchanged and on the remaining tensor product of  two copies of $V$ we have to define a map from $V\otimes V$ to $V\{n\}$. We define it on the basis vectors as: 
$$X^i \otimes X^j \mapsto X^{i+j},\quad i,j=0,\ldots,n,$$ 
where we assume $X^i=0$, for $i>n$.\\
\indent In the second case, the added edge belongs to one connected component, say $p$-th one, and the remaining components will remain unchanged. In this case, we have to define a degree preserving map from $V^{\otimes l}\{ni\}$ to $V^{\otimes l}\{n(i+1)\}$. We define the map as the identity on those $l-1$ copies of $V$ which correspond to the connected components which are not affected by adding the new edge, while on the remaining copy of $V$ (the one which corresponds to the $p$-th connected component), we have to define a linear map from $V$ to $V\{n\}$. We define it on the basis vectors by $1 \mapsto X^{n}$, and $X^i \mapsto 0$, for $i=1,\ldots,n$. \\
\begin{remark}
We could also define per-edge differentials in the same manner as in Section 3, by using the alternative desription of $V$ as the quotient of the ring of polynomials. Namely if the added edge connects two components, than the differential is given on the corresponding tensor product of two copies of $V$ by multiplication of polynomials (with setting $X_q$ to be $X_p$), and in the other case is given by multiplication by $X_p^n$.
\end{remark}
\indent These per-edge maps obviously make the cube commutative and after defining signs in the completely same way as before, and adding over the columns (as in figure \ref{sl2}), we obtain the differential $d$ of our chain complex $D_n(G)$. Again, as before we have that the homology groups $H^i(D_n(G))$ do not depend on the ordering of edges of the graph $G$, and since the differential is degree preserving, we have:

\begin{theorem}
The homology groups of the chain complex $D_n(G)$ are graph invariants. The graded Euler characteristic of the complex $D_n(G)$ is equal to the $n$-th specialization $Q_n(G)$ of the dichromatic polynomial of the graph $G$.
\end{theorem}

\end{document}